\documentclass[11pt]{article} 
 
\usepackage{latexsym} 
\usepackage{amsfonts} 
\usepackage{amsmath} 
\usepackage{amsthm} 
\usepackage{mathrsfs} 
\usepackage{dsfont} 
\usepackage{bbold} 
\usepackage[english]{babel} 
\usepackage{caption} 
\usepackage{epsfig} 
\usepackage{float} 
\usepackage{subfigure} 
\usepackage{psfrag} 
\usepackage{graphicx} 
\usepackage{epsfig} 
\usepackage{amsbsy} 
\usepackage{hyperref}
 
\DeclareMathOperator{\Val}{\matV}

\newtheorem{theorem}{Theorem} 
\newtheorem*{prop*}{Theorem} 
 
\newtheorem{coro}[theorem]{Corollary} 
 
\newtheorem{lemma}[theorem]{Lemma} 
 
\newtheorem{rmk}[theorem]{Remark}

\newtheorem{hyp}{Hypothesis}

\newcommand{\zerarcounters}{\setcounter{equation}{0}\setcounter{theorem}{0}} 
 
\newcommand{\ZZZ}{\mathds{Z}} 
 
\newcommand{\NNN}{\mathds{N}} 
 
\newcommand{\RRR}{\mathds{R}} 
\newcommand{\TTT}{\mathds{T}}

\newcommand{\calG}{{\mathcal G}}

\newcommand{\TT}{{\mathcal T}}

\newcommand{\gotn}{{\mathfrak n}}

\newcommand{\gotB}{{\mathfrak B}}

\newcommand{\gotT}{{\mathfrak T}}

\newcommand{\matV}{{\mathscr V}}

\newcommand{\ol}{\overline} 
 
\newcommand{\Fullbox}{{\rule{2.0mm}{2.0mm}}} 
 
\newcommand{\EP}{\hfill\Fullbox\vspace{0.2cm}} 
\newcommand{\prova}{\noindent{\it Proof. }} 
\newcommand{\io}{\infty} 
\newcommand{\e}{\varepsilon} 
\newcommand{\al}{\alpha}

\newcommand{\x}{\xi}

\newcommand{\g}{\gamma}

\newcommand{\oo}{\boldsymbol{\omega}} 
\newcommand{\mm}{\boldsymbol{\mu}}

\newcommand{\nn}{\boldsymbol{\nu}} 

\newcommand{\pps}{\boldsymbol{\psi}} 
\newcommand{\vzero}{\boldsymbol{0}}

\newcommand{\ii}{{\rm i}}

\def\tilde#1{\widetilde{#1}}
 
\def\ins#1#2#3{\vbox to0pt{\kern-#2 \hbox{\kern#1 #3}\vss}\nointerlineskip} 
 
\begin{document}
 
\title{\bf Convergent series for quasi-periodically forced strongly dissipative systems}
 
\author 
{\bf Livia Corsi$^{1}$, Roberto Feola$^{2}$ and Guido Gentile$^{3}$
\vspace{2mm} 
\\ \small
$^{1}$ Dipartimento di Matematica, Universit\`a di
Napoli ``Federico II'', Napoli, I-80126, Italy
\\ \small 
$^{2}$ Dipartimento di Matematica, Universit\`a di Roma ``La Sapienza", Roma, I-00185, Italy
\\ \small 
$^{3}$ Dipartimento di Matematica, Universit\`a di Roma Tre, Roma, I-00146, Italy
\\ \small 
E-mail:  livia.corsi@unina.it, feola@mat.uniroma1.it, gentile@mat.uniroma3.it}

\date{} 
 
\maketitle 
 
\begin{abstract} 
We study the ordinary differential equation
$\e \ddot x+\dot x + \e \, g(x) = \e f(\oo t)$,
with $f$ and $g$ analytic and $f$ quasi-periodic in $t$
with frequency vector $\oo\in\RRR^{d}$.
We show that if there exists $c_{0}\in\RRR$ such that
$g(c_{0})$ equals the average of $f$ and the first non-zero derivative
of $g$ at $c_{0}$ is of odd order $\gotn$, then, for $\e$ small enough
and under very mild Diophantine conditions on $\oo$, 
there exists a quasi-periodic solution close to $c_{0}$, 
with the same frequency vector as $f$. In particular
if $f$ is a trigonometric polynomial the Diophantine condition
on $\oo$ can be completely removed.
This extends results previously available in the literature for $\gotn=1$.
We also point out that, if  $\gotn=1$ and the first derivative of $g$ at $c_{0}$
is positive, then the quasi-periodic solution is locally unique and attractive.
 \end{abstract} 
  
\zerarcounters 
\section{Introduction} 
\label{sec:1} 

Consider the ordinary differential equation
\begin{equation}\label{eq:1.1}
\e \ddot x + \dot x + \e \, g(x) = \e \, f(\oo t) ,
\end{equation}
where $x\in\RRR$, $\e\in\RRR$, $\oo \in \RRR^{d}$, with $d\in\NNN$,
and the functions $g:\RRR\to\RRR$ and $f:\TTT^{d}\to\RRR$ are real analytic.
 From a physical point of view, for $\e>0$ the equation describes a one-dimensional system
with mechanical force $g$, subject to a quasi-periodic forcing $f$
with frequency vector $\oo$ and in the presence of friction --- with $1/\e$ being the
damping coefficient. The parameter $\e$ plays the r\^ole of a perturbation parameter.
Equations like (\ref{eq:1.1}) describe electric circuits
which are of interest in electronic engineering and theory of circuits;
we refer to \cite{TPJ,ADH,MCT,DM,BDGM}
for physical motivations and more details.

Without loss of generality we can (and shall) assume $\oo\cdot\nn\neq0$
$\forall\nn\in\ZZZ^{d}_{*}:=\ZZZ^{d}\setminus\{\vzero\}$.
If not, $f$ can be expressed as a quasi-periodic function with frequency vector
$\oo'\in\RRR^{d'}$, $d'<d$, with rationally independent components.
Let us denote by $\Sigma_{\xi}$ the strip of $\TTT^{d}$ of width $\xi$
and by $\Delta(c_{0},\rho)$ the disk of center $c_{0}$ and radius $\rho$
in the complex plane. By the assumptions on $f$ and $g$, for any $c_{0}\in\RRR$
there exist $\xi_{0}>0$ and $\rho_{0}>0$ such that
$f$ is analytic in $\Sigma_{\xi_{0}}$ and $g$ is analytic in $\Delta(c_{0},\rho_{0})$.
Then for all $\xi<\xi_{0}$ and all $\rho<\rho_{0}$ one has
\begin{subequations} \label{eq:1.2}
\begin{align}
f(\pps) & = \sum_{\nn\in\ZZZ^{d}} {\rm e}^{\ii\nn\cdot\pps} f_{\nn} , \qquad 
|f_{\nn}| \le \Phi \, {\rm e}^{-\x |\nn|} ,
\label{eq:1.2a} \\
g(x) & =\sum_{p=0}^{\io} a_{p} x^{p} , \qquad 
a_{p} := \frac{1}{p!} \frac{{\rm d}^{p}g}{{\rm d}x^{p}}(c_{0}) , \qquad
|a_{p}| \le \Gamma \, \rho^{\,p} ,
\label{eq:1.2b}
\end{align}
\end{subequations}
where $\Phi$ is the maximum of $f(\pps)$ for $\pps\in\Sigma_{\xi}$
and $\Gamma$ is the maximum of $g(x)$ for $x\in\Delta(c_{0},\rho)$.
Of course both $\rho_{0}$ and $\Gamma$ depend on $c_{0}$.

Let us define
\begin{equation} \label{eq:1.3}
\al_{n}(\oo) := \min \big\{ |\oo\cdot\nn| : 0<|\nn| \le 2^{n}  \big\} ,
\qquad \gotB(\oo) := \sum_{n=0}^{\io} \frac{1}{2^{n}} \log \frac{1}{\al_{n}(\oo)} ,
\end{equation}
and set
\begin{equation} \label{eq:1.4}
\beta_{n}(\oo) := \min \big\{ |\oo\cdot\nn| : 0<|\nn| \le 2^{n} , \; f_{\nn} \neq 0 \big\} ,
\qquad \e_{n}(\oo) := \frac{1}{2^{n}} \log \frac{1}{\beta_{n}(\oo)} .
\end{equation}

We shall require a \textit{non-resonance} condition on $\oo$ and
a \textit{non-degeneracy} condition relating $f$ to $g$, as follows.

\begin{hyp} \label{hyp:1}
Assume $\displaystyle{\lim_{n\to\io}\e_{n}(\oo)=0}$.
\end{hyp}

\begin{hyp} \label{hyp:2}
There exists a zero $c_{0}$ of odd order $\gotn$ of the equation $g(c_{0})=f_{\vzero}$.
\end{hyp}

Hypothesis \ref{hyp:1} is automatically satisfied either for $d=1$ (periodic case) or
for $d>1$ and $f$ a trigonometric polynomial. If $\oo$ is a Bryuno vector \cite{B},
that is if $\gotB(\oo)<\io$, then the sequence $\{\e_{n}(\oo)\}$ is summable,
so, for such $\oo$, Hypothesis \ref{hyp:1} is satisfied for any analytic function $f$.
Hypothesis \ref{hyp:2} implies that
\begin{equation} \label{eq:1.5}
\frac{{\rm d}^{k}g}{{\rm d}x^{k}}(c_{0}) = 0 \hbox{ for } k=1,\ldots,\gotn-1, \qquad
a:=\frac{1}{\gotn!} \frac{{\rm d}^{\gotn}g}{{\rm d}x^{\gotn}}(c_{0}) \neq 0 .
\end{equation}

The following result is proved in Section \ref{sec:2}.

\begin{theorem} \label{thm:1}
Assume Hypotheses \ref{hyp:1} and \ref{hyp:2} with $\gotn=1$.
There exists $\e_{0}>0$ such that for all $|\e|\le\e_{0}$ there is
at least one quasi-periodic solution $x_{0}(t)=c_{0}+u(\oo t,\e)$
to (\ref{eq:1.1}), with $u(\pps,t)=O(\e)$ analytic in $\pps$ and $C^{\io}$ in $\e$.
\end{theorem}

The periodic case is much easier and has been
studied in detail \cite{GBD1,GBD2,CCD},
with a  thorough characterization
of the analyticity properties for $\e$ in the complex plane.
One could study the analyticity properties of
the quasi-periodic solution for $\e$ in the complex plane
also in the quasi-periodic case --- for instance
this has been done in \cite{GBD1,GBD2,CCD,CFG}.
Here we prefer to focus on real $\e$, both for simplicity and because
it represents the interesting case from a physical point of view:
$\e>0$ small corresponds to a system with large damping coefficient $\g=1/\e$.

The existence of a quasi-periodic solution with the same frequency vector
as the forcing was proved in \cite{GBD2}
under the the stronger assumption that $\oo$ is a Bryuno vector.
The condition on $\oo$ has been weakened into Hypothesis \ref{hyp:1}
in \cite{CCD}, where the solution was also proved to be jointly analytic in $\pps$
and in $\e$, for $\e$ in a suitable domain of the complex plane
with the boundary tangent to the imaginary axis at the origin.
In Section \ref{sec:2} we give a different proof with respect to \cite{CCD},
where a method based on a fixed point argument is used,
by showing that, by introducing an auxiliary parameter $\mu$,
it is possible to write the solution as a convergent power series in $\mu$.
However we want to stress that we borrow from \cite{CCD} the
idea not to estimate all small divisors independently of $\e$:
we refer to Remark \ref{rmk:2.9} for a more precise comparison.

Both the method used here --- and in \cite{GBD1,GBD2,G1,G3} --- and the method
of \cite{CCD} assure the uniqueness of the solution
only in a suitable space of smooth functions,
so in principle we can neither exclude the existence of other quasi-periodic solutions
nor conclude that the solution we construct is attractive.
However, under a slightly stronger non-degeneracy condition we can
obtain more information. Indeed the following result
holds --- the proof can be found in \cite{G3}.

\begin{theorem} \label{thm:2}
Consider (\ref{eq:1.1}) with $\e>0$.
Assume Hypothesis \ref{hyp:2} with $\gotn=1$ and $a>0$.
If there is a quasi-periodic solution to (\ref{eq:1.1}) of the form
$x_{0}(t)=c_{0}+O(\e)$, then it is a local attractor.
\end{theorem}

By existence of local attractor we mean that there is
a simply connected open set containing the solution such that all trajectories
starting inside that set tend to the solution as time goes to infinity.
In particular this yields that the quasi-periodic solution is locally unique.
Combining Theorems \ref{thm:1} and \ref{thm:2} we deduce the following result.

\begin{theorem} \label{thm:3}
Consider (\ref{eq:1.1}) with $\e>0$ and $f$ a trigonometric polynomial. If 
Hypothesis \ref{hyp:2} holds with $\gotn=1$ and $a>0$, then for
any $\oo$ there is a quasi-periodic local attractor with the same
frequency vector as $f$. 
\end{theorem}

On physical grounds we could expect a result of this kind to hold for
any analytic function $f$.  However Theorem \ref{thm:1} requires
some Diophantine condition on $\oo$ --- however mild it may be.
It would be interesting to see whether the Diophantine condition on $\oo$
can be removed completely for $f$ analytic, as in the case
of trigonometric polynomials. Another interesting question is whether
results of the kind of Theorems \ref{thm:1} to \ref{thm:3} could
be obtained when $\gotn \ge 3$.

As far as the second question is concerned, the first remark in order is that
the assumption that the zero is of odd order cannot be removed:
if there is a zero $c_{0}$ of even order, then there is
no quasi-periodic solution to (\ref{eq:1.1}) reducing to $c_{0}$ as $\e\to0$ \cite{G1}.
For odd $\gotn$ we shall prove a result analogous to Theorem \ref{thm:1}.

\begin{theorem} \label{thm:4}
Assume Hypotheses \ref{hyp:1} and \ref{hyp:2} with $\gotn \ge 3$.
There exist $\e_{0}>0$ such that for all $|\e|\le\e_{0}$ there is
at least one quasi-periodic solution $x_{0}(t)=c_{0}+u(\oo t,\e)$
to (\ref{eq:1.1}), with $u(\pps,t)=O(\e)$ analytic in $\pps$ and $C^{\io}$ in $\e$.
\end{theorem}

If we require for $\oo$ to be a Bryuno vector, then the existence of
an analytic quasi-periodic solution of the form $c_{0}+O(\e)$
for all $\gotn$ odd follows from \cite{G1,G2}. The proof of Theorem \ref{thm:4}
follows the same lines of the proof of Theorem \ref{thm:1},
after a first step of perturbation theory in order to modify the linear operator,
and with a more careful use of the irrationality of
the frequency vector $\oo$. In Section \ref{sec:3} we shall give discuss
the case of trigonometric polynomials, and we shall see in Appendix \ref{app:b}
how to generalise the proof to any analytic $f$.
It would be worth investigating whether
some analogues to Theorems \ref{thm:2} and \ref{thm:3}
could be obtained for $\gotn \ge 3$.

\zerarcounters 
\section{Proof of Theorem \ref{thm:1}}
\label{sec:2} 

Here we assume both Hyothesis \ref{hyp:1} and Hypothesis \ref{hyp:2},
with $\gotn=1$ and define $a$ as in (\ref{eq:1.5}).
Let us rewrite (\ref{eq:1.1}) as
\begin{equation} \label{eq:2.1}
\e \ddot x + \dot x + \e \, a \, x + \e \, G(x) = \e \,  \tilde f(\oo t) ,
\end{equation}
where
\begin{equation} \label{eq:2.2}
G(x) := g(x)-g(c_{0})-a\, x = \sum_{p=2}^{\io} a_{p} x^{p} , \qquad
\tilde f(\pps) := f(\pps) - f_{\vzero} = \sum_{\nn\in\ZZZ^{d}_{*}} {\rm e}^{\ii\nn\cdot\pps} f_{\nn} ,
\end{equation}
and introduce the auxiliary parameter $\mu$ by modifying (\ref{eq:2.1}) into
\begin{equation} \label{eq:2.3}
\e \ddot x + \dot x + \e \, a \, x + \mu \, \e \, G(x) = \mu \, \e \, \widetilde f(\oo t) .
\end{equation}
Then we look for a quasi-periodic solution to (\ref{eq:2.3}) of the form
\begin{equation} \label{eq:2.4}
x(t,\e,\mu) = c_{0} + u(\oo t, \e,\mu) , \qquad u(\pps,\e,\mu) = \sum_{k=1}^{\io} \sum_{\nn\in\ZZZ^{d}} 
\mu^{k} {\rm e}^{\ii\nn\cdot\pps} u^{(k)}_{\nn}(\e) .
\end{equation}
We shall show that there exists $\mu_{0}>0$ such that 
there exists a solution of the form (\ref{eq:2.4}), analytic in $\mu$ for $|\mu|<\mu_{0}$.
Since the original equation is recovered when $\mu=1$ we need $\mu_{0}>1$. This will be
obtained by showing that the coefficients $u^{(k)}_{\nn}(\e)$ are bounded as
$|u^{(k)}_{\nn}(\e)| \le AB^{k}{\rm e}^{-\xi'|\nn|} |\e|^{\al k}$, for suitable positive constants $A,B,\xi',\al$.

By inserting (\ref{eq:2.4}) into (\ref{eq:2.3}) we obtain a recursive definition
for the coefficients $u^{(k)}_{\nn}(\e)$. By defining
\begin{equation} \label{eq:2.5}
D(\e,s) := - \e s^{2} + \ii s + \e \, a ,
\end{equation}
one has, formally,
\begin{subequations} \label{eq:2.6}
\begin{align}
D(\e, \oo\cdot\nn) \, u^{(1)}_{\nn} (\e) & = \e \, f_{\nn} 
\label{eq:2.6a} \\
D(\e, \oo\cdot\nn) \, u^{(k)}_{\nn} (\e) & = - \e \, \sum_{p=2}^{\io} a_{p} \!\!\!\!\!\!
\sum_{\substack{ k_{1},\ldots,k_{p} \ge 1 \\ k_{1}+\ldots+k_{p}=k-1}}
\sum_{\substack{\nn_{1},\ldots,\nn_{p} \in\ZZZ^{d} \\ \nn_{1}+\ldots+\nn_{p}=\nn}} \!\!\!\!\!\!
u^{(k_{1})}_{\nn_{1}} (\e) \ldots u^{(k_{p})}_{\nn_{p}} (\e) , \qquad k \ge 2 ,
\label{eq:2.6b}
\end{align}
\end{subequations}
for $\nn\neq0$, and
\begin{equation} \label{eq:2.7}
a \, u^{(k)}_{\vzero}(\e)  = - 
 \sum_{p=2}^{\io} a_{p} \!\!\!\!\!\!
 \sum_{\substack{ k_{1},\ldots,k_{p} \ge 1 \\ k_{1}+\ldots+k_{p}=k-1}}
\sum_{\substack{\nn_{1},\ldots,\nn_{p} \in\ZZZ^{d} \\ \nn_{1}+\ldots+\nn_{p}=\vzero}} \!\!\!\!\!\!
u^{(k_{1})}_{\nn_{1}} (\e) \ldots u^{(k_{p})}_{\nn_{p}} (\e) , \qquad k \ge 1 .
\end{equation}
Here and henceforth the sums over the empty set are meant as zero.

\begin{rmk} \label{rmk:2.1}
\emph{
For $k=1$ (\ref{eq:2.7}) yields $u^{(1)}_{\vzero}=0$.
For $k=2$ one has $u^{(2)}_{\nn}=0$ $\forall\nn\in\ZZZ^{d}$.
}
\end{rmk}

By iterating the definition one obtains an explicit expression for the coeffiecients
$u^{(k)}_{\nn}$, which can be represented in terms of trees.

A \textit{rooted tree} $\theta$ is a graph with no cycle,
such that all the lines are oriented toward a unique
point (\textit{root}) which has only one incident line (root line).
All the points in $\theta$ except the root are called \textit{nodes}.
The orientation of the lines in $\theta$ induces a partial ordering 
relation ($\preceq$) between the nodes. Given two nodes $v$ and $w$,
we shall write $w \prec v$ every time $v$ is along the path
(of lines) which connects $w$ to the root; we shall write $w\prec \ell$ if
$w\preceq v$, where  $v$ is the unique node that $\ell$ exits.
For any node $v$ denote by $p_{v}$ the number of lines entering $v$.

Given a rooted tree $\theta$ we denote by $N(\theta)$ the set of nodes,
by $E(\theta)$ the set of \textit{end nodes}, i.e. nodes $v$ with $p_{v}=0$,
by $V(\theta)$ the set of \textit{internal nodes}, i.e. nodes $v$ with $p_{v}\ge 1$,
and by $L(\theta)$ the set of lines; by definition $N(\theta)=E(\theta) \amalg V(\theta)$.

We associate with each end node $v\in E(\theta)$
a \textit{mode} label $\nn_{v}\in\ZZZ^{d}_{*}$
and with each internal node an \textit{degree} label $d_{v}\in\{0,1\}$.
With each line $\ell\in L(\theta)$ we associate
a \textit{momentum} $\nn_{\ell} \in \ZZZ^{d}$ with the constraint
\begin{equation} \label{eq:2.8}
\nn_{\ell}=\sum_{\substack{w\in E(\theta) \\ w \prec \ell}} \nn_{w} .
\end{equation}
We add the two following further constraints:
(1) $p_{v}\ge 2$ $\forall v\in V(\theta)$ and
(2) if $d_{v}=0$ then the line $\ell$ exiting $v$ has $\nn_{\ell}=\vzero$.
We shall write $V(\theta)=V_{0}(\theta) \amalg V_{1}(\theta)$,
where $V_{0}(\theta):=\{ v\in V(\theta): d_{v}=0\}$.
For any discrete set $A$ we denote by $|A|$ its cardinality.
Define the \textit{degree} and the \textit{order} of $\theta$ as
$d(\theta):=|E(\theta)|+|V_{1}(\theta)|$ and $k(\theta):=|N(\theta)|$, respectively.

We call \textit{equivalent} two labelled rooted trees which can be transformed into
each other by continuously deforming the lines in such a way that
they do not cross each other. In the following we shall consider only
inequivalent labelled rooted trees, and we shall call them call trees \textit{tout court},
for simplicity.

We associate with each node $v\in N(\theta)$  a \textit{node factor}
\begin{equation} \label{2.9}
F_{v} := \begin{cases}
- \e^{d_{v}}  \, a_{p_{v}} , & v \in V(\theta) , \\
\e \, f_{\nn_{v}} , & v \in E(\theta) ,
\end{cases}
\end{equation}
and with each line $\ell\in L(\theta)$ a \textit{propagator}
\begin{equation} \label{eq:2.10}
\calG_{\ell} := \begin{cases}
1/D(\e, \oo\cdot\nn_{\ell}) , & \nn_{\ell} \neq \vzero , \\
1/a , & \nn_{\ell}=\vzero .
\end{cases}
\end{equation}
Then, by defining
\begin{equation} \label{eq:2.11}
\Val(\theta,\e) := \Biggl( \prod_{v\in N(\theta)} F_{v} \Biggr) \Biggl( \prod_{\ell\in L(\theta)} \calG_{\ell} \Biggr)
\end{equation}
one has
\begin{equation} \label{eq:2.12}
u^{(k)}_{\nn}(\e) = \sum_{\theta\in \TT_{k,\nn}} \Val(\theta,\e) , \quad \nn \in\ZZZ^{d} ,
\end{equation}
where $\TT_{k,\nn}$ is the set of trees of order $k$
and momentum $\nn$ associated with the root line.

\begin{lemma} \label{lem:2.2}
One has $|D(\e,s)| \ge \max\{|a\e|,|s|\}$ for $\e$ small enough and all $s\in\RRR$.
\end{lemma}

\prova
One has $|D(\e,s)| \ge |{\rm Im}\,D(\e,s)|$ and $|D(\e,s)|\ge |D(\e,0)|$ for $\e$ small enough.
\EP

\begin{lemma} \label{lem:2.3}
For any tree $\theta$ one has $|E(\theta)|\ge |V(\theta)|+1$.
\end{lemma}

\prova
By induction on the order of the tree.
\EP

\begin{rmk} \label{rmk:2.4}
\emph{
Equality $|E(\theta)|=|V(\theta)|+1$ holds when $|N(\theta)|=2^{p}+1$, with $p\ge 1$, 
and $p_{v}=2$ for all $v\in V(\theta)$.
}
\end{rmk}

\begin{coro} \label{coro:2.5}
For any tree $\theta$ one has $\displaystyle{|E(\theta)|\ge \frac{1}{2} \left( k(\theta)+1 \right)}$.
\end{coro}

\begin{lemma} \label{lem:2.6}
For any $k\ge 1$, any $\nn\in\ZZZ^{d}$ and any tree $\theta\in\TT_{k,\nn}$ one has
\begin{equation} \nonumber
\left| \Val(\theta,\e) \right| \le A B^{k}  |\e|^{(k+1)/2} 
\prod_{v\in E(\theta)} {\rm e}^{-3\xi |\nn_{v}|/4}
\end{equation}
where $\xi$ is as in (\ref{eq:1.2a}),
with $A=1$ and $B$ a positive constant depending on $\Phi$, $\Gamma$ and $\rho$.
\end{lemma}

\prova
One bounds (\ref{eq:2.11}) as
\begin{equation} \nonumber
\left| \Val(\theta,\e) \right| \le |\e|^{d(\theta)} \Biggl( \prod_{v\in V(\theta)} | a_{p_{v}} | \Biggr) 
\Biggl( \prod_{v \in E(\theta)}  \frac{| f_{\nn_{v}} | }{ |\oo\cdot\nn_{v}|} \Biggr)
\Biggl( \prod_{v\in V_{0}(\theta)} \frac{1}{|a|} \Biggr)
\Biggl( \prod_{v\in V_{1}(\theta)} \frac{1}{|a\e|} \Biggr) ,
\end{equation}
where we have used the bound $|D(\e,s)|\ge |s|$ for the propagators
of the lines exiting the end nodes and the bound $|D(\e,s)|\ge |a\e|$
for the propagators of the lines exiting the nodes in $V_{1}(\theta)$.
For each end node we bound $f_{\nn_{v}}$ as in (\ref{eq:1.2}): 
then we extract a factor ${\rm e}^{-3\xi|\nn_{v}|/4}$
and use Hypothesis \ref{hyp:1} to bound ${\rm e}^{-\xi|\nn|/4}|\oo\cdot\nn|^{-1} \le C_{0}$,
for a suitable constant $C_{0}$, for all $\nn$ such that $f_{\nn} \neq 0$.
Moreover, by Corollary \ref{coro:2.5},
\begin{equation} \nonumber
d(\theta) -|V_{1}(\theta)| = |E(\theta)| \ge \frac{k(\theta)+1}{2} ,
\end{equation}
so that we obtain
\begin{equation} \nonumber
\left| \Val(\theta,\e) \right| \le \Gamma^{|V(\theta)|} \rho^{|N(\theta)|}
(C_{0}\Phi)^{|E(\theta)|} a^{-|V(\theta)|} {\rm e}^{-3\xi |\nn|/4} 
|\e|^{(k(\theta)+1)/2} .
\end{equation}
Therefore, by bounding $\max\{|E(\theta)|,|V(\theta)|\} \le k(\theta)$,
the assertion follows.
\EP

\begin{lemma} \label{lem:2.7}
For any $k\ge 1$ and $\nn\in\ZZZ^{d}$ one has
\begin{equation} \nonumber
\left| u^{(k)}_{\nn}(\e) \right| \le A\, C^{k} {\rm e}^{-\xi |\nn|/2} |\e|^{(k+1)/2} ,
\end{equation}
where $\xi$ is as in (\ref{eq:1.2}),
with $A=1$ and $C$ a positive constant depending on $\Phi$, $\Gamma$, $\xi$ and $\rho$.
\end{lemma}

\prova
The coefficients $u^{(k)}_{\nn}$ are defined by (\ref{eq:2.12}):
we have to use the bounds of Lemma \ref{lem:2.6}
and sum over all trees in $\TT_{k,\nn}$. The sum over 
the Fourier labels $\{\nn_{v}\}_{v\in E(\theta)}$
is performed thanks to the factors ${\rm e}^{-3\xi |\nn_{v}|/4}$ 
associated with the end nodes that we have not used to control the 
denominators $|\oo\cdot\nn_{v}|$ --- see the proof
of Lemma \ref{lem:2.6} --- and produces an overall factor
$C_{1}^{|E(\theta)|}{\rm e}^{-\xi |\nn|/2}$, for some positive constant $C_{1}$. 
The sum over the other labels produces a factor $C_{2}^{|N(\theta)|}$,
with $C_{2}$ a suitable positive constant.
Then the assertion follows by taking $C=BC_{1}C_{2}$.
\EP

\begin{rmk} \label{rmk:2.8}
\emph{
The main idea in the proof is to bound in a different way the propagators,
depending on whether or not the lines exit end nodes.
Eventually the propagators of the lines exiting the end nodes have 
a ``gain" factor $\e$ with respect to the propagators of the other lines:
together with the fact that each internal node has at least
two entering lines --- so that the number of ``bad'' propagators turns
out to be less than the number of ``good" propagators ---, this
leads to bound the product of the propagators 
of any tree of order $k$ proportionally to $|\e|^{-k/2}$.
Note that a similar feature has been exploited in \cite{GGG}
in a rather different context, i.e. the problem of synchronisation
in chaotic systems. As in that case --- and as in \cite{CCD} --- no small divisor problem 
arises: of course this makes easier to study the convergence of the series.
}
\end{rmk}

\begin{rmk} \label{rmk:2.9}
\emph{
The crucial property described in Remark \ref{rmk:2.8},
which allows to require only Hypothesis \ref{hyp:1} on $\oo$,
has been already pointed out and used in \cite{CCD}:
in our proof we simply adapted that idea to our formalism.
Smoothness in $\e$ at $\e=0$ is not discussed in \cite{CCD},
but very likely could be derived also with the method used therein.
}
\end{rmk}

The function (\ref{eq:2.4}), with the coefficients given by (\ref{eq:2.12}),
solves (\ref{eq:2.3}) order by order. Since the series converges uniformly,
then it is also a solution \textit{tout court} of (\ref{eq:2.3}) --- and hence
of (\ref{eq:2.1}) for $\mu=1$.  Anayticity in $\pps\in\Sigma_{\xi'}$
for any $\xi'<\xi/2$ follows from the bound on the Fourier coefficients
given by Lemma \ref{lem:2.7}.
To prove smoothness in $\e$ one can reason as follows.
Each value $\Val(\theta,\e)$ is a polynomial in $\e$ with coefficients
depending on $\e$ through the propagators. If one compute the
$n$-th derivative of $\Val(\theta,\e)$ with respect to $\e$,
one can bound it in a different way depending on whether
one has or not $(k(\theta)+1)/2 \le n$. If $(k(\theta)+1)/2 \le n$,
then all the propagators and their derivatives
are bounded by using the inequality $|D(\e,s)|\ge |s|$ of Lemma \ref{lem:2.2}.
If $(k(\theta)+1)/2>n$ one can reason as done above to arrive at
the bounds in Lemma \ref{lem:2.7}: one obtains the same bounds,
with a coefficient $A$ depending on $n$
and with a power of $\e$ decreased by $n$, so that the sum over
$k(\theta)$ can still be performed.

\zerarcounters 
\section{Proof of Theorem \ref{thm:4}}
\label{sec:3} 

Assume Hypothesis \ref{hyp:1} and Hypothesis \ref{hyp:2} with $\gotn\ge 3$.
We look for a solution $x(t)$ to (\ref{eq:1.1}) of the form
\begin{equation} \label{eq:3.1}
x(t) = c_{0} + \e \, x_{1}(t) + \xi(t) ,
\end{equation}
where $\e\, x_{1}(t)$ it the solution to the first-order truncation of (\ref{eq:1.1}), i.e.
\begin{equation} \nonumber
\e \ddot x_{1} + \dot x_{1} = \tilde f(\oo t) ,
\end{equation}
with $\tilde f$ as in (\ref{eq:2.2}). An easy computation 
gives $x_{1}(t) = \zeta + u^{[1]}(\oo t,\e)$, where
\begin{equation} \label{eq:3.2}
u^{[1]}(\pps,\e) := \sum_{\nn\in\ZZZ^{d}_{*}} {\rm e}^{\ii\nn\cdot\pps}
u^{[1]}_{\nn}(\e) , \qquad u^{[1]}_{\nn}(\e):=
\frac{f_{\nn}}{\ii\oo\cdot\nn (1+\ii \e \oo\cdot\nn)} ,
\end{equation}
and $\zeta$ is a real parameter that will be fixed later on. 
Note that Hypothesis \ref{hyp:1} guarantees that the
function (\ref{eq:3.2}) is well-defined for any analytic $f$.
Therefore the problem is reduced to finding a zero-average
quasi-periodic solution $\xi(t)$ to  the equation
\begin{equation} \label{eq:3.3}
\e \ddot \xi + \dot \xi + \e \, \tilde G(\e \, x_{1}(t)+\xi) = 0 , 
\qquad \tilde G(x) := \sum_{p=\gotn}^{\io} a_{p} x^{p} ,
\end{equation}
which can be rewritten as
\begin{equation} \label{eq:3.4}
\e \ddot \xi + \dot \xi + b \, \e^{\gotn} \, \xi + 
\mu \, \e \, \widehat G( \mu \e x_{1}(t), \xi) = 0 ,
\end{equation}
where $\mu=1$ and
\begin{subequations} \label{eq:3.5}
\begin{align}
b & := \sum_{p=\gotn}^{\io} p \, a_{p} \, \e^{p-\gotn} 
\bigl[ \bigl( x_{1}(t) \bigr)^{p-1} \bigr]_{\vzero} ,
\label{eq:3.5a} \\
\qquad \widehat G(x,\xi) & := 
\sum_{p=\gotn}^{\io} a_{p} \sum_{s=0}^{p}
\left( \begin{matrix} p \\ s \end{matrix} \right)
\xi^{s} \left(  x^{p-s} - \delta_{s,1} \left[ x^{p-s} \right]_{\vzero} \right) ,
\label{eq:3.5b}
\end{align}
\end{subequations}
with $[\cdot]_{\vzero}$ denoting --- here and henceforth --- the average on $\TTT^{d}$.

\begin{rmk} \label{rmk:3.1}
\emph{
By setting $b_{0}:=
\gotn \, a_{\gotn}\,  [ (x_{1}(t) )^{\gotn-1} ]_{\vzero}$,
one has $b_{0}\neq0$, because $a_{\gotn}=a \neq0$ and $\gotn$ is odd,
and hence $b=b_{0}(1+O(\e))$ does not vanish for $\e$ small enough.
}
\end{rmk}

As in Section \ref{sec:2} we first ignore the constraint $\mu=1$ and treat it as a parameter:
we shall look for a solution which can be written as a power series in $\mu$, 
with coefficients which still admit a tree expansion.

Let us assume here that $f$ is a trigonometric polynomial of degree $N$,
i.e. that $f_{\nn}=0$ for all $\nn\in\ZZZ^{d}$ such that $|\nn|>N$. 
In such a case it is more convenient to redefine
$\Phi=\max\{|f_{\nn}| : |\nn| \le N\}$. We shall see in Appendix \ref{app:b}
how to extend the proof to the case of $f$ analytic. Define
\begin{equation} \label{eq:3.6}
\al = \min \{ |\oo\cdot \nn| : 0< |\nn| \le (\gotn+1)\, N \} .
\end{equation}
One has $\al>0$ by the assumption of irrationality on $\oo$.

With respect to Section \ref{sec:2} we modify the tree expansion as follows.
Rooted trees and the sets $N(\theta)$, $E(\theta)$, $V(\theta)$ and $L(\theta)$
are defined as previously. 
If $p_{v}$ denotes the number of lines entering $v\in V(\theta)$ we impose
the constraint $p_{v}\ge \gotn$. 
We associate with each end node $v$ a mode label $\nn_{v}\in\ZZZ^{d}$
and with each line a momentum $\nn_{\ell}\in\ZZZ^{d}$,
still satisfying (\ref{eq:2.8}) and with
the further constraint that $\nn_{\ell}\neq\vzero$ if $\ell$ exits an internal node.
We split $E(\theta)=E_{0}(\theta)\amalg E_{1}(\theta)$, with
$E_{0}(\theta):=\{v\in E(\theta) : \nn_{v}=\vzero\}$.
The order of $\theta$ is still defined as $k(\theta)=|N(\theta)|$.
Let us define also $L_{0}(\theta):=\{\ell \in L(\theta) : |\oo\cdot\nn_{\ell}|<\al/2\}$
and $V_{0}(\theta):=\{v\in V(\theta) : \hbox{the line $\ell$ exiting $v$ belongs
to $L_{0}(\theta)$}\}$, and set $V_{1}(\theta):=V(\theta)\setminus V_{0}(\theta)$.

We call \textit{excluded} a node $v$ such that
$p_{v}-1$ lines entering $v$ exit end nodes, and the other line entering $v$
exits an internal node and has the same momentum as the line exiting $v$.
Let $\gotT_{k,\nn}$ be the set of inequivalent labelled rooted trees,
which do not contain any excluded nodes,
of order $k$ and momentum $\nn$ associated with the root line.
In the following we shall call simply trees  the elements of $\gotT_{k,\nn}$.

We associate with each node $v\in N(\theta)$ a \textit{node factor}
\begin{equation} \label{3.7}
F_{v} := \begin{cases}
- \e \, a_{p_{v}} , & v \in V(\theta) , \\
\e \, f_{\nn} , & v \in E_{1}(\theta) \\ 
\e \, \zeta , & v \in E_{0}(\theta) ,
\end{cases}
\end{equation}
where $\zeta\in\RRR$ is the parameter introduced before (\ref{eq:3.2}),
and with each line $\ell\in L(\theta)$ a \textit{propagator}
\begin{equation} \label{eq:3.8}
\calG_{\ell} := \begin{cases}
\calG_{E}(\e, \oo\cdot\nn_{\ell}) , & 
\nn_{\ell} \neq \vzero  \hbox{ and $\ell$ exits an end node} , \\
\calG_{V}(\e, \oo\cdot\nn_{\ell}) , & 
\nn_{\ell} \neq \vzero \hbox{ and $\ell$ exits an internal node} , \\
1 , & \nn_{\ell}=\vzero , \\
\end{cases}
\end{equation}
with
\begin{equation} \label{eq:3.9}
\calG_{E}(\e,s) := \frac{1}{\ii s (1+\ii \e s)} , \qquad
\calG_{V}(\e,s) := \frac{1}{\ii s (1+\ii \e s) + b\,\e^{\gotn}} ,
\end{equation}
where $b\in\RRR_{+}$ is defined in (\ref{eq:3.5a}) --- and hence is a function of $\zeta$.

Setting
\begin{equation} \label{eq:3.10}
\Val(\theta,\e) := \Biggl( \prod_{v\in N(\theta)} F_{v} \Biggr) 
\Biggl( \prod_{\ell\in L(\theta)} \calG_{\ell} \Biggr)
\end{equation}
and
\begin{equation} \label{eq:3.11}
u^{[k]}_{\nn}(\e) := \!\!\!\! \sum_{\theta\in\gotT_{k,\nn}} \!\! 
\Val(\theta,\e) , \quad \nn \neq \vzero  ,
\qquad k \ge 2 ,
\end{equation}
we define (formally) the series
\begin{equation} \label{eq:3.12}
\ol{\xi}(\pps,\e,\mu) := 
\sum_{k=2}^{\io} \mu^{k} \sum_{\nn\in\ZZZ^{d}_{*}} {\rm e}^{\ii\nn\cdot\pps} u^{[k]}_{\nn}(\e) .
\end{equation}
and set $\ol{u}(\pps,\e,\mu)=\mu \e \bigl( \zeta + u^{[1]}(\pps,\e) \bigr) + \ol{\xi}(\pps,\e,\mu)$.

\begin{rmk} \label{rmk:3.2}
\emph{
The constraint $p_{v} \ge \gotn$ implies 
$u^{[k]}_{\nn}(\e)=0$ $\forall\nn\in\ZZZ^{d}_{*}$ and $2\le k \le \gotn$.
}
\end{rmk}

\begin{rmk} \label{rmk:3.3}
\emph{
The coefficients (\ref{eq:3.10}) depend on $\zeta$, which so far is still a free parameter.
}
\end{rmk}

The definition (\ref{eq:3.12}) is formal not only in the sense that it may fail to converge.
In fact the very definition of the coefficients $u^{[k]}_{\nn}(\e)$ involves
quantities --- the propagators --- for which we have not yet any estimate.
The latter problem is easily solved as follows. Define
\begin{equation} \label{eq:3.13}
D_{V}(\e,s) := \frac{1}{\calG_{V}(\e,s)} = \ii s \left( 1 + \ii \e s \right) + b\, \e^{\gotn} .
 \end{equation}
 %

\begin{lemma} \label{lem:3.4}
There exist $\e_{1}>0$ such that
$|D_{N}(\e,s)| \ge \max\{ |s|,|b\e^{\gotn}|  \}$ for all $s\in\RRR$ and $|\e|<\e_{1}$.
\end{lemma}

\prova
Reason as in the proof of Lemma \ref{lem:2.2}.
\EP

Thanks to Lemma \ref{lem:3.4} we deduce that the coefficients (\ref{eq:3.11})
are well defined for all $k\ge 2$ and all $\nn\in\ZZZ^{d}_{*}$.
Now we want to find conditions for the series (\ref{eq:3.12}) to converge.
We shall prove that, for any $\zeta\in\RRR$, under the assumption
that $\e$ is small enough, depending on $\zeta$,
the series (\ref{eq:3.12}) converges to a well-defined function
analytic in $\pps$ and $C^{\io}$ in $\e$,
with a radius of convergence $\mu_{0}>1$: 
this will allow us to take $\mu=1$ in (\ref{eq:3.4}).
Moreover we shall show that, for any fixed $\ol{\zeta}\in\RRR$,
the coefficients $u^{[k]}_{\nn}(\e)$ admit uniform bounds for $|\zeta|\le \ol{\zeta}$.

\begin{lemma} \label{lem:3.5}
For any tree $\theta$ one has $|E(\theta)| \ge (\gotn-1)\,|V(\theta)| + 1$.
\end{lemma}

\prova
The bound is proved by induction by using that $p_{v}\ge \gotn$ $\forall v\in V(\theta)$.
\EP

\begin{coro} \label{coro:3.6}
For any tree $\theta$ one has  $\gotn\, |E(\theta)| \ge (\gotn-1)\,k(\theta) +1$.
\end{coro}

\begin{lemma} \label{lem:3.7}
For any tree $\theta$ one has $\gotn\,|V_{0}(\theta)| \le E(\theta) - 2$.
\end{lemma}

The proof of Lemma \ref{lem:3.7} is in Appendix \ref{app:a}.

\begin{lemma} \label{lem:3.8}
For any $k\ge 1$ and $\nn\in\ZZZ^{d}_{*}$ and any tree $\theta\in\gotT_{k,\nn}$ one has
\begin{equation} \nonumber
\left| \Val(\theta,\e) \right| \le A\, B^{k}
|\zeta|^{|E_{0}(\theta)|} |b|^{-|V_{0}(\theta)|} 
|\e|^{1+\frac{\gotn-1}{\gotn^{2}} k} ,
\end{equation}
with $A=1$ and $B$ a positive constant depending on $\Phi$, $\Gamma$, $\rho$ and $\al$.
\end{lemma}

\prova
One bounds (\ref{eq:3.10}) as
\begin{equation} \nonumber
\left| \Val(\theta,\e) \right| \le |\e|^{k(\theta)}  |\zeta|^{|E_{0}(\theta)|}
\Biggl( \prod_{v\in V(\theta)} \!\!\!\! |a_{p_{v}} | \Biggr) 
\Biggl( \prod_{v \in E_{1}(\theta)}  \frac{| f_{\nn_{v}} | }{ |\oo\cdot\nn_{v}|} \Biggr)
\Biggl(  \prod_{v \in V_{0}(\theta)}   \frac{1}{| b\e^{\gotn}| }\Biggr) 
\Biggl( \prod_{v \in V_{1}(\theta)} \frac{2}{\al} \Biggr) ,
\end{equation}
where, by relying on Lemma \ref{lem:3.4},
we have used the bound $|D(\e,s)|\ge |s|$ for the propagators
of the lines exiting either the nodes $v \in E_{1}(\theta)$ or the nodes $v\in V_{1}(\theta)$
and the bound $|D(\e,s)|\ge |b\e^{\gotn}|$ for the propagators
of the lines exiting the nodes in $V_{0}(\theta)$. For each end node
we bound $|f_{\nn_{v}}|\,|\oo\cdot\nn_{v}|^{-1} \le 2\Phi/\alpha$. Then
\begin{equation} \nonumber
\left| \Val(\theta,\e) \right| \le |\e|^{k(\theta)-\gotn|V_{0}(\theta)|}
\Gamma^{|V(\theta)|} \rho^{|N(\theta)|}
|\zeta|^{|E_{0}(\theta)|} \Phi^{|E_{1}(\theta)|}
|b|^{-|V_{0}(\theta)|} (2/\al)^{|V_{1}(\theta)|+|E_{1}(\theta)|} ,
\end{equation}
where we can bound, by using
Lemma \ref{lem:3.7} and Corollary \ref{coro:3.6},
\begin{eqnarray}
k(\theta)-\gotn |V_{0}(\theta)| & \!\!\! = \!\!\! &
|E(\theta)|+|V(\theta)|-\gotn |V_{0}(\theta)|  \ge
|E(\theta)|-(\gotn -1) |V_{0}(\theta)| , \nonumber \\
& \!\!\! \ge \!\!\! &
 |E(\theta)|-(\gotn -1) \frac{E(\theta) - 2}{\gotn} =
2 -\frac{2}{\gotn} + \frac{|E(\theta)|}{\gotn} \ge
1 + \frac{\gotn-1}{\gotn^{2}} k(\theta) , \nonumber
\end{eqnarray}
so that, by using that $\max\{|V(\theta)|,|E(\theta)|\} \le |N(\theta)|=k(\theta)$,
the assertion  follows.
\EP

\begin{rmk} \label{rmk:3.9}
\emph{
The bounds in Lemma \ref{lem:3.8} depend on $\zeta$. However,
for any given $\ol{\zeta}>0$, there exists $\ol{b}>0$ such that
$|b|\ge \ol{b}$ (by compactness), and hence for $|\zeta|\le \ol{\zeta}$
we can obtain uniform bounds
\begin{equation} \nonumber
| \Val(\theta,\e)| \le A\, B^{k}
\zeta_{*}^{k} b_{*}^{-k/\gotn}|\e|^{1+\frac{\gotn-1}{\gotn^{2}} k} ,
\end{equation}
where $\zeta_{*}=\max\{1,\ol{\zeta}\}$ and $b_{*}=\min\{1,\ol{b}\}$.
}
\end{rmk}

\begin{lemma} \label{lem:3.10}
For any $k\ge 1$ and $\nn\in\ZZZ^{d}_{*}$ one has
\begin{equation} \nonumber
\left| u^{[k]}_{\nn}(\e) \right| \le A\, C^{k}
|\e|^{1+\frac{\gotn-1}{\gotn^{2}} k} ,
\end{equation}
with $A=1$ and $C$ a positive constant depending on $\Phi$, $\Gamma$, $\rho$,
$\zeta_{*}$, $b_{*}$, $N$ and $\al$.
\end{lemma}

\prova
Reason as in the proof of Lemma \ref{lem:2.7} and use Remark \ref{rmk:3.9}.
Now the sum over the mode labels can be bounded by $(2N+1)^{d k}$.
\EP

We have proved that the series (\ref{eq:3.12}) converges to
a well-defined a function for $|\mu|<\mu_{0}$, with $\mu_{0}>1$, provided
(1) $|\zeta| \le \ol{\zeta}$ for some $\ol{\zeta}>0$ and
(2) $\e$ is small enough.
Moreover, by construction, the function
is periodic and analytic in $\pps$ and $C^{\io}$ in both $\e$ and $\zeta$
(this can be seen as in Section \ref{sec:2}).
In the remaining part of this section we shall prove that 
$\zeta$ can be fixed in such a way that $|\zeta| \le \ol{\zeta}$
and the function $\ol{\xi}(\oo t,\e,\mu)$ solves the equation (\ref{eq:3.4})
for $|\mu|<\mu_{0}$ --- and hence in particular
$\ol{\xi}(\oo t,\e,1)$ solves the equation (\ref{eq:3.3}).

We can write (\ref{eq:3.4}) in Fourier space, if we expand (formally)
\begin{equation} \nonumber
\xi(t) = \sum_{\nn\in\ZZZ^{d}_{*}} {\rm e}^{\ii\nn\cdot \oo t} \xi_{\nn} ,
\end{equation}
so as to obtain
\begin{subequations} \label{eq:3.14}
\begin{align}
\left( \ii\oo\cdot\nn (1+\ii \e\oo\cdot\nn) + b\, \e^{\gotn} \right) \xi_{\nn}  & =
- \bigl[ \e \, \widehat G(\mu \, \e \, x_{1}(t),\xi ) \bigr]_{\nn} , \qquad \nn\neq\vzero ,
\label{eq:3.14a} \\
0 & = - \bigl[ \e \, \widehat G(\mu \, \e \, x_{1}(t),\xi ) \bigr]_{\vzero}  ,
\label{eq:3.14b}
\end{align}
\end{subequations}
where $[A]_{\nn}$ means that we expand the function $A$ in Fourier series in $\pps$
and keep the Fourier coefficient with label $\nn$. If we expand further
(again formally) $\xi_{\nn}$ as a Taylor series in $\mu$, by writing
\begin{equation} \nonumber
\xi_{\nn} = \sum_{k=2}^{\io} \mu^{k} \xi_{\nn}^{[k]} ,
\end{equation}
we can write (\ref{eq:3.14a}) order by order,
\begin{equation} \label{eq:3.15}
\left( \ii\oo\cdot\nn (1+\ii \e\oo\cdot\nn) + b\, \e^{\gotn} \right)
\xi^{[k]}_{\nn} = - \bigl[ \e \, \widehat G(\mu \, \e \, x_{1}(t),\xi ) 
\bigr]^{[k-1]}_{\nn} , \quad \nn \neq \vzero , \quad k \ge 2, 
\end{equation}
where $[A]^{[k]}_{\nn}$ means that we expand the function $[A]_{\nn}$
in powers of $\mu$ and keep the Taylor coefficient to order $k$.

\begin{lemma} \label{lem:3.11}
For any $\zeta\in\RRR$ the coefficients $\xi^{[k]}_{\nn} = u^{[k]}_{\nn}(\e)$
solve (\ref{eq:3.15}).
\end{lemma}

\prova
Expand the right hand side of (\ref{eq:3.15}) in powers
of $\e x_{1}$ and $\xi$, and write $\e\,x_{1}$ according
to (\ref{eq:3.1}), with the coefficients $u^{[1]}_{\nn}(\e)$ as in (\ref{eq:3.2}),
and $\xi$ according to (\ref{eq:3.12}), with the coefficients $u^{[k]}_{\nn}(\e)$ 
as in (\ref{eq:3.11}). Then (\ref{eq:3.15}) reduces to (\ref{eq:3.11}) itself.
\EP

Therefore, for any $\zeta\in\RRR$, the function (\ref{eq:3.12}) formally solves
(\ref{eq:3.14a}). If $|\zeta|\le \ol{\zeta}$ and $\e$ is small enough
then the series (\ref{eq:3.12}) converges uniformly and therefore solves (\ref{eq:3.14a}).
Moreover the function $\widehat G(\mu \, \e \, x_{1}(t),\xi )$ is well defined
and hence it makes sense to consider its average.
So we are left with the equation (\ref{eq:3.14b}): we shall show that
it is possible to fix $\zeta$ in such a way that such an equation is satisfied.

Consider the implicit function problem
\begin{equation} \label{eq:3.16}
F_{2}(\zeta,\e) : = \frac{1}{\e^{\gotn}} 
\bigl[ \e \, \widehat G(\mu \, \e \, x_{1}(t),\xi ) \bigr]_{\vzero} = 0 .
\end{equation}
If we are able to find a solution to (\ref{eq:3.16}) then (\ref{eq:3.14b}) is also satisfied.

\begin{lemma} \label{lem:3.12}
Let $\e$ be small enough.
There is $\ol{\zeta}>0$ such that there exists a unique real solution
$\zeta=\zeta(\e)$ to (\ref{eq:3.16}), with $|\zeta(\e)|<\ol{\zeta}$.
Moreover $\zeta(\e)$ is $C^{\io}$ in $\e$.
\end{lemma}

\prova
One has
\begin{equation} \nonumber
F_{2}(\zeta,\e) = \ol{F}_{2}(\zeta)  + F_{3}(\zeta,\e) , \qquad \ol{F}_{2} (\zeta) := 
\big[ \left( \zeta + u^{[1]}(\oo t,\e)  \right)^{\gotn} \big]_{\vzero}  ,
\end{equation}
where $F_{3}(\zeta,\e)$ is a function which goes to  zero when $\e$ goes to zero.
The function  $\ol{F}_{2}(\zeta)$ is a polynomial of order $\gotn$
in $\zeta$. The equation $\ol{F}_{2}(\zeta)=0$ admits
a unique real root $\zeta_{0}$ --- see Lemma 2.4
in \cite{G2}. Since ${\rm d} \ol{F}_{2}(\zeta)/{\rm d}\zeta=\gotn
\big[ \left( \zeta + u^{[1]}(\oo t,\e) \right)^{\gotn-1} \big]_{\vzero} = b_{0}/a$,
the root is simple. Therefore, by the implicit function theorem,
for $\e$ small enough there is a unique $\zeta(\e)$ such that
$\zeta(0)=\zeta_{0}$ and $F_{2}(\zeta(\e),\e)=0$.
\EP

In particular, as the proof of Lemma \ref{lem:3.12} shows,
one can take $\ol{\zeta}=2\zeta_{0}$, where $\zeta_{0}$ is the root
of $\ol{F}_{2}(\zeta)=0$. The proof of Theorem \ref{thm:4}
is complete, in the case of trigonometric polynomials $f$.

\zerarcounters 
\section{Comments}
\label{sec:4} 

By looking at the proofs of Theorems \ref{thm:1} and \ref{thm:4} given in
Sections \ref{sec:2} and \ref{sec:3}, respectively, we see that, in both cases,
Hypothesis \ref{hyp:1} ensures that the first order is well-defined.
Indeed we used Hypothesis \ref{hyp:1} to bound the propagators
of the lines $\ell$ exiting the end nodes as $|\oo\cdot\nn_{\ell}|^{-1}$,
i.e. to control $D(\e,\oo\cdot\nn)$ in (\ref{eq:2.6a}) and
the denominators $\ii\oo\cdot\nn(1+\ii\e\oo\cdot\nn)$ in (\ref{eq:3.3}).
Therefore we can rephrase the two theorems by saying that
the condition for a quasi-periodic to exist is the same condition
for the first order to be well-defined.

We note also that, under Hypothesis \ref{hyp:1}, the formal expansion
in powers of $\e$ is well defined to all orders. This can be easily checked,
for instance, by looking for a solution to (\ref{eq:1.1})
in the form of a formal power series in $\e$:
one finds a recursive definition for the coefficients of the series and writes down
a tree expansion for such coefficients. Then one shows that the coefficients are well
defined to all orders under the only Hypothesis \ref{hyp:1} --- see \cite{CG}
for more details in a similar case. In our case, the existence of the formal power series
is enough to conclude about the existence of a solution: Theorems \ref{thm:1}
and \ref{thm:4} imply that the conditions for the existence of a solution
are the the same conditions required for the existence of a formal solution,
i.e. a solution in the form of a formal power series in $\e$.
This is a quite non-general feature. Usually, one cannot infer that a solution
to an ordinary differential equation exists simply
from the  fact that a formal solution exists: a classical example
are elliptic lower-dimensional tori (see however \cite{S}
and references quoted therein for germs of vector fields
in a neighborhood of a fixed point with hyperbolic linear part).

 From the perspective outlined above
 it is not clear whether the theorems could be generalised
to any vector $\oo\in\RRR^{d}$ for any analytic forcing $f$:
for the first order or the formal expansion
to be well defined, some Diophantine condition on $\oo$ seems necessary.
On the other hand, for instance when $\gotn=1$
and $a>0$, from a physical point of view one expects for a local attractor to exist,
and it is not unlikely that a synchronisation phenomenon occurs. Of course
it could also happen that the conjugation exists but is not smooth
(think of Denjoy's theorem for diffeomorphisms of the circle \cite{D,A}): in that case
it would not not possible to construct it with the techniques used here.

We conclude with two technical comments.
\begin{enumerate}
\item The smallness assumptions of both Theorems \ref{thm:1} and \ref{thm:4}
can be weakened. Indeed, by looking at the proofs of the Theorems --- in
particular the bounds on $\Val(\theta,\e)$ given in the proofs of Lemmas
\ref{lem:2.6} and \ref{lem:3.8} ---, one sees that the parameter
that must be small is $\e\Phi$. Therefore
a large forcing is still allowed as far as $\e\Phi$ remains small.
\item A property like Lemma \ref{lem:3.7} --- or Lemma \ref{lem:b.5}
in the case of analytic forcing --- can be found to hold
also for the case $\gotn=1$ (for a suitably defined
set $V_{0}(\theta)$ --- or $L_{0}(\theta)$ in Appendix \ref{app:b}).
The argument of Section \ref{sec:2} shows that this is
not necessary to prove the existence of a quasi-periodic solution.
However, a property of this kind can be used to enlarge,
with respect to the results found in \cite{CCD},
the domain of analyticity for $\e$ in the complex plane;
see \cite{CFG} for results in that direction.
\end{enumerate}

\appendix

\zerarcounters 
\section{Proof of Lemma \ref{lem:3.7}}
\label{app:a} 

The proof is by induction on the order of the tree. First of all note that
$k(\theta) \ge \gotn+1$ by construction (see Remark \ref{rmk:3.2}).
If $k(\theta)=\gotn +1$, then the root line has momentum
$\nn=\nn_{1}+\ldots+\nn_{\gotn}$, where $\nn_{1},\ldots,\nn_{\gotn}$
are the mode labels of the $\gotn$ end nodes of $\theta$, so that  $|\nn|\le \gotn N$
and hence $|\oo\cdot\nn| \ge \al > \al/2$. Therefore $|V_{0}(\theta)|=0$
in such a case and the bound holds.

If $k(\theta)\ge \gotn+2$, call $\ell_{0}$ the root line of $\theta$
and $v_{0}$ the node which $\ell_{0}$ exits. Let
$r$ be the number of end nodes whose exiting lines
enter $v_{0}$ and set $s=p_{v_{0}}-r\,$: there will be
$s$ trees $\theta_{1},\ldots,\theta_{s}$ such that the
respective root lines $\ell_{1},\ldots,\ell_{s}$ enter $v_{0}$.
Note that $|E(\theta)|=|E(\theta_{1})|+\ldots+|E(\theta_{s})|+r$.
If $s=0$ the bound holds: indeed if $k(\theta)=\gotn+2$, then
$|E(\theta)|=\gotn+1$ and $|V_{0}(\theta)|=0$,
because the momentum $\nn$ of the root line is such that
$|\nn| \le (\gotn+1)N$ and hence $|\oo\cdot\nn|\ge\al$, while
if $k(\theta)\ge \gotn+3$, then $|E(\theta)|\ge \gotn+2$ and $|V_{0}(\theta)|\le 1$.
Therefore in the following we assume $s\ge 1$.

If $\ell_{0} \notin L_{0}(\theta)$, then $|V_{0}(\theta)|=|V_{0}(\theta_{1})|+\ldots+
|V_{0}(\theta_{s})|$, so that, by the inductive hypothesis,
\begin{equation} \nonumber
|V_{0}(\theta)| \le \sum_{k=1}^{s} \frac{E(\theta_{k})-2}{\gotn}
\le \frac{E(\theta)- r -  2 s }{\gotn}
\le \frac{E(\theta) -  2}{\gotn} ,
\end{equation}
and the bound follows.

If $\ell_{0}\in L_{0}(\theta)$, then $|V_{0}(\theta)|=1+|V_{0}(\theta_{1})|+\ldots+
|V_{0}(\theta_{s})|$ and, again by the inductive hypothesis,
\begin{equation} \nonumber
|V_{0}(\theta)| \le 1 + \sum_{k=1}^{s} \frac{E(\theta_{k})-2}{\gotn}
\le \frac{E(\theta) -  2}{\gotn} + \Biggl[ 1 -
\frac{r  + 2\,(s-1) }{\gotn} \Biggr] ,
\end{equation}
If $s+r\ge \gotn+1$, then $r + 2 \, (s-1) \ge \gotn + (s-1)\ge \gotn$.
If $s+r=\gotn$ and $s\ge 2$, then $r + 2\,(s-1) \ge 
\gotn+(s-2) \ge \gotn$. Thus in both cases
the last term in square brackets  is non-positive and the bound follows.

If $s+r=\gotn$ and $s=1$, then the line $\ell_{1}$ must be
in $L_{1}(\theta)$. This can be seen by reductio ad absurdum.
Suppose that $\ell_{1}\in L_{0}(\theta)$. Then
$|\oo\cdot\nn_{\ell_{1}}|< \al/2$. Moreover
$|\oo\cdot\nn_{\ell_{0}}|<\al/2$ because $\ell_{0}\in L_{0}(\theta)$ by hypothesis.
On the other hand one has $\nn_{\ell_{0}}=\nn_{\ell_{1}}
+\nn_{1}+\ldots+\nn_{r}$, where $r=\gotn-1$ and $\nn_{1},\ldots,\nn_{r}$
are the mode labels of the $r$ end nodes whose exiting lines enter $v_{0}$.
Therefore, if we use that $|\nn_{1}+\ldots+\nn_{r}| \le (\gotn-1)N$
and $\nn_{1}+\ldots+\nn_{r} \neq \vzero$ (otherwise $v_{0}$
would be an excluded node), we obtain
\begin{equation} \nonumber
\al > |\oo\cdot\nn_{\ell_{0}}| + |\oo\cdot\nn_{\ell_{1}}| \ge
|\oo\cdot (\nn_{\ell_{0}}-\nn_{\ell_{1}})| =
|\oo\cdot (\nn_{1}+\ldots+\nn_{r})| \ge \al ,
\end{equation}
so arriving at a contradiction.  Let $v_{1}$ be the node
which $\ell_{1}$ exits: there will be $r'$ end nodes whose exiting lines
enter $v_{1}$ and $s'$ trees $\theta_{1}',\ldots,\theta_{s'}'$ whose root lines 
$\ell_{1}',\ldots,\ell_{s}'$ enter $v_{1}$. One has
$|E(\theta)|=|E(\theta_{1}')| + \ldots + |E(\theta_{s'}')| + r + r'$
and $|V_{0}(\theta)|=1 + |V_{0}(\theta_{1}')| + \ldots + |V_{0}(\theta_{s'}')|$,
where $s=1$, $r=\gotn-1$ and $r'+s'\ge \gotn$.
By the inductive hypothesis one has
\begin{equation} \nonumber
|V_{0}(\theta)| \le
1 + \sum_{k'=1}^{s'} \frac{E(\theta_{k'}')-2}{\gotn} \le
\frac{E(\theta)  -  2}{\gotn} + \Biggl[ 1 -
\frac{r +r'+2s'-2}{\gotn} \Biggr] ,
\end{equation}
where $r +r'+2s'-2 \ge 2\gotn+s'-3\ge \gotn +(\gotn-3) \ge \gotn$,
so that the last term in square bracket is non-positive.
Therefore the bound follows once more.

\zerarcounters 
\section{Proof of Theorem \ref{thm:4} for analytic forcing}
\label{app:b} 

In the analytic case the trees are constructed as in Section \ref{sec:3}:
in particular the definition of the coefficients (\ref{eq:3.11}) of
the series (\ref{eq:3.12}) is the same. The only difference is how
to bound the values of the trees in (\ref{eq:3.11}).

First of all we need some notations. We shall not introduce the sets
$V_{0}(\theta)$ and $V_{1}(\theta)$ of Section \ref{sec:3}. Instead,
we shall proceed as follows. For any node $v\in V(\theta)$ define
$E(\theta,v):=\{w\in E(\theta) : \hbox{ the line exiting $w$ enters $v$}\}$,
$r_{v}:=|E(\theta,v)|$, $s_{v}:=p_{v}-r_{v}$ and
\begin{equation} \nonumber
\mm_{v}:=\sum_{w\in E(\theta,v)} \nn_{w}  ,
\qquad \mu_{v}:=|\mm_{v}| .
\end{equation}
Set $V_{0}(\theta):=\{v\in V(\theta) : s_{v}=0\}$
and $V_{1}(\theta):=\{v\in V(\theta) : s_{v}=1\}$.
If $v\in V_{0}(\theta)$ we call $\ell_{v}$ the line exiting $v$, while
if $v\in V_{1}(\theta)$ we call $\ell_{v}$ the line exiting $v$ and $\ell_{v}'$
the line entering $v$ which does not exits an end node.

\begin{rmk} \label{rmk:b.1}
\emph{
By Hypothesis \ref{hyp:1} there exists $C_{0}>0$ such that
$C_{0}|\oo\cdot\nn| \ge {\rm e}^{-\xi|\nn|/8}$ $\forall \nn\in\ZZZ^{d}_{*}$.
}
\end{rmk}
\begin{lemma} \label{lem:b.2}
If $v\in V_{0}(\theta)$ one has $C_{0}|\oo\cdot\nn_{\ell_{v}}|
\ge {\rm e}^{-\xi \mu_{v}/8}$.
\end{lemma}

\prova
For $v\in V_{0}(\theta)$ one has $\nn_{\ell_{v}}=\mm_{v}$,
so that the bound follows from Remark \ref{rmk:b.1}.
\EP

\begin{lemma} \label{lem:b.3}
If $v\in V_{1}(\theta)$ one has $C_{0}\max\{|\oo\cdot\nn_{\ell_{v}}|,
|\oo\cdot\nn_{\ell_{v}'}|\} \ge {\rm e}^{-\xi \mu_{v}/8}/2$.
\end{lemma}

\prova
By contradiction: if the bound does not hold then
\begin{equation}
{\rm e}^{-\xi \mu_{v}/8}  > C_{0} |\oo\cdot\nn_{\ell_{v}}| + 
C_{0} |\oo\cdot\nn_{\ell_{v}'}|
\ge C_{0} |\oo\cdot(\nn_{\ell_{v}}-\nn_{\ell_{v}'})|= C_{0}
| \oo\cdot \mm_{v} | \ge {\rm e}^{-\xi \mu_{v}/8} ,
\end{equation}
where we have used that $\mm_{v} \neq \vzero$, 
otherwise $v$ would be an excluded node.
\EP

Define $L_{1}(\theta,v):=\{\ell_{v}\}$ for $v\in V_{0}(\theta)$ and
$L_{1}(\theta,v):=\{\ell\in\{\ell_{v},\ell_{v}'\}:
C_{0} |\oo\cdot\nn_{\ell}| \ge {\rm e}^{-\xi \mu_{v}/8}/2\}$ for $v\in V_{1}(\theta)$.
By Lemmas \ref{lem:b.2} and \ref{lem:b.3} one has
$L_{1}(\theta,v)\neq\emptyset$ for all $v\in V_{0}(\theta)\cup V_{1}(\theta)$.
Set also $L_{1}(\theta):=\{ \ell \in L(\theta) : \exists v \in
V_{0}(\theta)\cup V_{1}(\theta)\hbox{ such that } \ell \in L_{1}(\theta,v)\}$,
$L_{\rm int}(\theta):=\{\ell\in L(\theta) : \ell \hbox{ exits
a node } v\in V(\theta)\}$ and $L_{0}(\theta) := L_{\rm int}(\theta)\setminus
L_{1}(\theta)$.

\begin{lemma} \label{lem:b.4}
For any tree $\theta$ one has $\gotn \, |L_{0}(\theta)| \le |E(\theta)| - 2$.
\end{lemma}

\prova
One proceeds by induction on $V(\theta)$.
If $|V(\theta)|=1$ then $V(\theta)=V_{0}(\theta)$ and hence
$|L_{0}(\theta)|=0$, while $|E(\theta)|-2>0$, so that the bound holds.
If $|V(\theta)|\ge 2$ the root line $\ell_{0}$ of $\theta$ exits a
node $v_{0}\in V(\theta)$ with $s_{v_{0}}+r_{v_{0}}\ge \gotn$
and $s_{v_{0}}\ge 1$. Call $\theta_{1},\ldots,\theta_{s_{v_{0}}}$
the trees whose respective root lines $\ell_{1},\ldots,\ell_{s_{v_{0}}}$
enter $v_{0}$: one has $|E(\theta)|=|E(\theta_{1})|+\ldots+|E(\theta_{s_{v_{0}}})|
+r_{v_{0}}$. If $\ell_{0}\notin L_{0}(\theta)$ then
$|L_{0}(\theta)|=|L_{0}(\theta_{1})|+\ldots+|L_{0}(\theta_{s_{v_{0}}})|$
and the bound follows from the inductive hypothesis.

If $\ell_{0}\in L_{0}(\theta)$ then one has
$|L_{0}(\theta)|=1 + |L_{0}(\theta_{1})|+\ldots+|L_{0}(\theta_{s_{v_{0}}})|$,
so that,  again by the inductive hypothesis,
\begin{equation} \nonumber
|L_{0}(\theta) \le \frac{|E(\theta)|-2}{\gotn} +
\Biggl[ 1 - \frac{r_{v_{0}} + 2\, (s_{v_{0}}-1)}{\gotn} \Biggl] ,
\end{equation}
so that, if either $r_{v_{0}}+s_{v_{0}}\ge \gotn+1$ or
$r_{v_{0}}+s_{v_{0}}=\gotn$ and $s_{v_{0}}\ge 2$, the bound follows.

If $r_{v_{0}}+s_{v_{0}}=\gotn$ and $s_{v_{0}}=1$, then $v_{0}\in V_{1}(\theta)$ 
and, since $C_{0}|\oo\cdot\nn_{\ell_{0}}|<{\rm e}^{-\xi\mu_{v_{0}}/8}/2$
(because $\ell\in L_{0}(\theta)$), then $C_{0}|\oo\cdot\nn_{\ell_{1}}| \ge
{\rm e}^{-\xi\mu_{v_{0}}/8}/2$ by Lemma \ref{lem:b.3}.
Therefore $\ell_{1}\notin L_{0}(\theta)$. If $v_{1}$ is the line which $\ell_{1}$
exits, call $\theta_{1}',\ldots,\theta_{s_{v_{1}}}'$ the trees whose root lines
enter $v_{1}$: one has
$|L_{0}(\theta)|=1+|L(\theta_{1}')|+\ldots+|L_{0}(\theta_{s_{v_{1}}}')|$
and hence, by the inductive hypothesis,
\begin{equation} \nonumber
|L_{0}(\theta)| \le 1 + \frac{|E(\theta)|-r_{v_{0}}-r_{v_{1}} -2s_{v_{1}}}{\gotn}
\le \frac{|E(\theta)|-2}{\gotn} +
\Biggl[ 1 - \frac{r_{v_{0}}+r_{v_{1}} +2\, (s_{v_{1}}-1)}{\gotn} \Biggr] ,
\end{equation}
where $r_{v_{0}}+r_{v_{1}}+2s_{v_{1}}-2 \ge \gotn$,
so that the bound follows in this case too.
\EP

\begin{lemma} \label{lem:b.5}
For any $k\ge 1$ and $\nn\in\ZZZ^{d}_{*}$ and any tree $\theta\in\gotT_{k,\nn}$
there are positive constants $A$ and $B$ such that
\begin{equation} \nonumber
\left| \Val(\theta,\e) \right| \le A\, B^{k} 
|\zeta|^{|E_{0}(\theta)|} |b|^{-|L_{0}(\theta)|}
|\e|^{1+\frac{\gotn-1}{\gotn^{2}} k}
\prod_{v\in E_{1}(\theta)} {\rm e}^{-5\xi |\nn_{v}|/8} ,
\end{equation}
where $\xi$ is as in (\ref{eq:1.2}),
with $A=1$ and the constant $B$ depending on $\Phi$, $\Gamma$ and $\rho$.
\end{lemma}

\prova
One bounds (\ref{eq:3.10}) as
\begin{equation} \nonumber
\left| \Val(\theta,\e) \right| \le |\e|^{k(\theta)}  |\zeta|^{|E_{0}(\theta)|}
\Biggl( \prod_{v\in V(\theta)} \!\!\!\! |a_{p_{v}} | \Biggr) 
\Biggl( \prod_{v \in E_{1}(\theta)}  | f_{\nn_{v}} | \Biggr)
\Biggl(  \prod_{\ell \in L(\theta)} \calG_{\ell} \Biggr) .
\end{equation}
We deal with the propagators as follows. If $\ell$ exits
a node $v\in V_{0}(\theta)$, then we have
\begin{equation} \nonumber
\left| \calG_{\ell} \right|
\!\!\!\!\! \prod_{w\in E_{1}(\theta,v)} \!\!\!
|f_{\nn_{w}} | \, |\calG_{\ell_{w}}| \le 
\frac{1}{|\oo\cdot\nn_{\ell}|}
\prod_{w\in E_{1}(\theta,v)} \frac{|f_{\nn_{w}}|}{|\oo\cdot\nn_{w}|} \le
C_{0}(\Phi C_{0})^{|E_{1}(\theta,v)|} \!\!\!\!\!\!
\prod_{w\in E_{1}(\theta,v)} \!\!\!\!\!\! {\rm e}^{-3\xi|\nn_{w}|/4} , 
\end{equation}
where $\ell_{w}$ is the line exiting the end node $w$ and
we have defined $E_{1}(\theta,v) :=\{w\in E(\theta,v) : \nn_{w} \neq \vzero\}$.
For the lines in $L_{1}(\theta)$ which do not exit nodes $v\in V_{0}(\theta)$
we distinguish three cases: given a node $v \in V_{1}(\theta)$
and denoting by $v'$ the node $\ell_{v}'$ exits,
(1) if one has $\ell_{v}\in L_{1}(\theta,v)$ and
either $\ell_{v}' \notin L_{1}(\theta,v)$ or $\ell_{v}' \in L_{1}(\theta,v')$,
we proceed as for the nodes $v\in V_{0}(\theta)$, so as to obtain
\begin{equation} \nonumber
\left| \calG_{\ell_{v}} \right| 
\!\!\!\!\! \prod_{w\in E_{1}(\theta,v)} \!\!\!
|f_{\nn_{w}} | \, |\calG_{\ell_{w}}| \le 
C_{0}(\Phi C_{0})^{|E_{1}(\theta,v)|} \!\!\!\!\!\!
\prod_{w\in E_{1}(\theta,v)} \!\!\!\!\!\! {\rm e}^{-3\xi|\nn_{w}|/4} ;
\end{equation}
(2) if $L_{1}(\theta,v)=\{\ell_{v}'\}$ and $\ell_{v'}\notin L_{1}(\theta,v')$, we bound
\begin{equation} \nonumber
\left| \calG_{\ell_{v}'} \right| 
\!\!\!\!\! \prod_{w\in E_{1}(\theta,v)} \!\!\!
|f_{\nn_{w}} | \, |\calG_{\ell_{w}}| \le 
C_{0}(\Phi C_{0})^{|E_{1}(\theta,v)|} \!\!\!\!\!\!
\prod_{w\in E_{1}(\theta,v)} \!\!\!\!\!\! {\rm e}^{-3\xi|\nn_{w}|/4} ;
\end{equation}
(3) if both lines $\ell_{v},\ell_{v}'$ belong to $L_{1}(\theta,v)$
and $\ell_{v}'\notin L_{1}(\theta,v')$, we bound
\begin{equation} \nonumber
\left| \calG_{\ell_{v}} \calG_{\ell_{v}'}\right|
\prod_{w\in E_{1}(\theta,v)} \frac{f_{\nn_{w}}}{|\oo\cdot\nn_{w}|} \le
C_{0}^{2}(\Phi C_{0})^{|E_{1}(\theta,v)|} \!\!\!\!\!\!
\prod_{w\in E_{1}(\theta,v)} \!\!\!\!\!\! {\rm e}^{-5\xi|\nn_{w}|/8} .
\end{equation}
For all the other propagators we bound
(1) $|\calG_{\ell}| \le 1$ if $\ell$ exits an end node $v$ with $\nn_{v}=\vzero$,
(2) $|\calG_{\ell}| \le |\oo\cdot\nn_{\ell}|^{-1}$ if $\ell$ exits an end node
$v$ with $\nn_{v}\neq\vzero$ and has not been already used in
the bounds above for the lines $\ell\in L_{1}(\theta)$, and
(3) $|\calG_{\ell}| \le |b\e^{\gotn}|^{-1}$ if $\ell \in L_{0}(\theta)$.
Then we obtain
\begin{equation} \nonumber
\left| \Val(\theta,\e) \right| \le |\e|^{k(\theta)-\gotn|L_{0}(\theta)|}
\Gamma^{|V(\theta)|} \rho^{|N(\theta)|}
|\zeta|^{|E_{0}(\theta)|} C_{0}^{|L_{1}(\theta)|}
(C_{0}\Phi)^{|E_{1}(\theta)|} |b|^{-|L_{0}(\theta)|} 
{\rm e}^{-5\xi |\nn|/8} ,
\end{equation}
where we can bound, by using Corollary \ref{coro:3.6} and
Lemma \ref{lem:b.5},
\begin{equation}
k(\theta)-\gotn |L_{0}(\theta)| =
|E(\theta)|+|V(\theta)|-\gotn |L_{0}(\theta)|  \ge
|E(\theta)|-(\gotn -1) |L_{0}(\theta)| \ge
1 + \frac{\gotn-1}{\gotn^{2}} k(\theta) , \nonumber
\end{equation}
so that the assertion  follows.
\EP

Fix $\ol{\zeta}$ and $\ol{b}$, and define $\zeta_{*}$ and $b_{*}$ as
in Section \ref{sec:3}.

\begin{lemma} \label{lem:b.6}
For any $k\ge 1$ and $\nn\in\ZZZ^{d}$ there are
positive constants $A$ and $C$ such that
\begin{equation} \nonumber
\left| u^{(k)}_{\nn}(\e) \right| \le A\, C^{k} {\rm e}^{-\xi |\nn|/2}
|\e|^{1+\frac{\gotn-1}{\gotn^{2}} k} ,
\end{equation}
where $\xi$ is as in (\ref{eq:1.2}), with $A=1$ and the constant $C$
depending on $\Phi$, $\Gamma$, $\rho$, $\xi$, $\zeta_{*}$ and $b_{*}$.
\end{lemma}

\prova
Reason as in the proof of Lemma \ref{lem:2.7}.
\EP

 From this point onward the proof proceeds as in the case of
trigonometric polynomial, so we skip the details.


\end{document}